\documentclass[12pt]{amsart}
\usepackage[english]{babel}
\usepackage[utf8]{inputenc}
\usepackage{amssymb,amsmath,amsfonts,amsthm,amscd,latexsym,graphicx,epsfig,colordvi}

\begin{document}

\title[On radicals of Novikov algebras]{On radicals of Novikov algebras}
\author{A.S.Panasenko}
\address{Alexander Sergeevich Panasenko
\newline\hphantom{iii} Sobolev Institute of Mathematics,
\newline\hphantom{iii} pr. Koptyuga, 4,
\newline\hphantom{iii} 630090, Novosibirsk, Russia.
\newline\hphantom{iii} Novosibirsk State University,
\newline\hphantom{iii} Universitetskiy pr., 1,
\newline\hphantom{iii} 630090, Novosibirsk, Russia
}
\email{a.panasenko@g.nsu.ru}%
\thanks{\sc Panasenko, A.S.,
On radicals of Novikov algebras}
\thanks{\copyright \ 2023 Panasenko A.S}

\maketitle {\small
\begin{quote}
\noindent{\sc Abstract. } We show that in a prime nonassociative Novikov algebra every nonzero ideal is non-associative. We prove that Baer (and Andrunakievich) radical and the largest left quasiregular ideal coincide in finite dimensional Novikov algebras over a field of characteristic 0 or algebraically closed field of odd characteristic. We show non-existence of right quasiregular radical in finite dimensional Novikov algebras. \medskip

 \end{quote}
}

\section{Introduction}

The structure theory of any variety of algebras is strongly related to the theory of radicals that exist in this variety. For example, in alternative algebras there are many useful radicals well known from the associative case: the Jacobson, Baer, Köthe, Levitsky, Andrunakievich radicals. In Jordan algebras the situation is slightly different: the quasiregular radical does not have such a good characterization as in the associative case, and the existence of the Baer radical has been an open question for more than half a century. However, in Jordan algebras there are a lot of other radicals (Andrunakievich, Köthe, Levitsky), as well as the McCrimmon radical, which is specific for the Jordan case.

In recent years, an active study of the structure and combinatorial theories of Novikov algebras has resumed. In the paper \cite{Bokut} L.~A.~Bokut, Y.~Chen and Z.~Zhang proved an analogue of the Poincaré-Birkhoff-Witt theorem for Novikov algebras and constructed the Gröbner-Shirshov basis for the free Novikov algebra.

In the article \cite{ShestZhang} I.~P.~Shestakov and Z.~Zhang proved that the following three properties are equivalent for the Novikov algebra $A$: $A$ is solvable, $A$ is right-nilpotent, and $A^2$ is nilpotent. In addition, in this paper it is shown that the minimal ideal of the Novikov algebra is either simple or has zero multiplication.

In \cite{UZ1}, U.~U.~Umirbaev and V.~N.~Zhelyabin proved an analog of the Bergman-Isaacs theorem for Novikov algebras: if $G$ is a finite abelian group and the zero component in a $G$-graded Novikov algebra $A$ is solvable, then the algebra $A$ is also solvable (if the characteristic of the field is zero or does not divide $|G|$).

Every left-symmetric algebra (and hence every Novikov algebra) is Lie-admissible. The paper \cite{UZ2} is devoted to the study of Lie-solvable Novikov algebras. For example, in this work K.~M.~Tulenbaev, U.~U.~Umirbaev and V.~N.~Zhelyabin showed that in a Lie-solvable Novikov algebra the commutator ideal is right-nilpotent (if the characteristic of the field is not equal to 2). The paper \cite{UZ3} is devoted to the study of simple Lie-solvable left-symmetric algebras.

In the paper \cite{P2022} the author shows that every ideal in a prime Novikov algebra is a prime Novikov algebra. A similar statement was proved for semiprime algebras, which made it possible to assert the existence of a Baer radical (in the sense of Kurosh) in the variety of Novikov algebras.

In this paper, radicals of Novikov algebras are considered. First, we refine the result from \cite{P2022} by proving that in a prime non-associative Novikov algebra, non-zero ideals are non-associative. This result will be used to characterize the Andrunakievich radical.

In many varieties of algebras, all popular and useful radicals coincide with the largest nilpotent ideal in the finite-dimensional case. It is well known that in finite-dimensional Novikov algebras, nilpotency is not a property that is stable under extensions. But solvability is. We prove that the Andrunakievich and Baer radicals coincide with each other and with the largest solvable ideal in the finite-dimensional Novikov algebra. Moreover, in the case of a field of characteristic 0 or an algebraically closed field of odd characteristic, the left quasiregular radical also coincides with the Baer radical.

Throughout this paper, all algebras are considered over a field $F$ of characteristic not 2. Definitions 2--7 are recalled from \cite{Zhevl}.

If $A$ is an algebra over a field, then for every elements $x,y,z\in A$ we use the notation $(x,y,z)=(xy)z-x(yz)$ for the \textbf{\textit{associator}} of elements $x,y,z$.

\textbf{Definition 1} (see \cite{Osborn}). An algebra $A$ over a field $F$ is called a \textbf{\textit{Novikov algebra}} if the following identities hold for all $x,y,z\in A$:
\[(x,y,z)=(y,x,z),\]
\[(xy)z=(xz)y.\]

\section{Semiprime Novikov algebras}

Let $A$ be a Novikov algebra. It is well known that the product of ideals in $A$ is an ideal. If $I$ is an ideal of the algebra $A$, then by $T(I)$ we denote the subspace in $I$ generated by all elements of the form $(ij)j$, where $i,j\in I$. Denote by $P(I)$ the subspace of $I$ generated by all elements of the form $(aj)j$, where $j\in I, a\in A$.\medskip

\textbf{Lemma 1.} \textit{Let $I$ be an ideal of the Novikov algebra $A$. Then the subspaces $P(I)$ and $T(I)$ are ideals in the algebra $A$, $P(I)=(AI)I$, $T(I)=I^2I$, and \[P (I)^2\subset T(I)\subset P(I).\]}

\textbf{Proof.} Let $i,j\in I$ and $a\in A$. Then $(a(i+j))(i+j), (ai)i, (aj)j\in P(I)$. Hence
\[2(ai)j=(ai)j+(aj)i=(a(i+j))(i+j)-(ai)i-(aj)j\in P(I).\]
Thus, $(ai)j\in P(I)$ and $P(I)=(AI)I$. Similarly, $T(I)=I^2I$. In particular, $T(I)$ and $P(I)$ are ideals of $A$.

It remains to note that if $i_1\in I, a_1\in A$, then
\[((aj)j)((a_1j_1)j_1)=((a((aj_1)j_1))j)j\in T(I).\]
So, $P(I)^2\subset T(I)$. \qed

\medskip

\textbf{Definition 2.} Recall that an algebra over a field is called \textbf{\textit{semiprime}} if every nonzero ideal $I$ has nontrivial multiplication, i.e. $I^2=0$ implies that $I=0$. An algebra is \textbf{\textit{prime}} if $IJ\neq 0$ for all nonzero ideals $I,J$. \medskip

\textbf{Lemma 2.} \textit{Let $I$ be a nonzero ideal of the semiprime Novikov algebra $A$. Then $T(I)\neq 0$.}

\textbf{Proof.} Assume that $T(I)=0$. Then, by Lemma 1, $P(I)^2=0$, whence, by semiprimality, $P(I)=0$, i.e.
\[(ai)j=0\]
for every $a\in A$, $i,j\in I$ by Lemma 1. In particular, $(I^2)^2=0$, so $I^2=0$ and $I= 0$. \qed

\medskip

\textbf{Definition 3.} If $A$ is an algebra over a field, then its \textbf{\textit{nucleus}} is the following subset:
\[N(A)=\{n\in A\mid (n,A,A)=(A,n,A)=(A,A,n)=0\}.\]
In every algebra, the nucleus is a subalgebra, see \cite{Zhevl}. In a Novikov algebra, the nucleus is an ideal, see \cite{P2022}. Moreover, in \cite{P2022} it is proved that the nucleus in a prime nonassociative Novikov algebra is equal to zero. \medskip

One can check that any Novikov algebra satisfies the following identities (see, for example, \cite{UZ2}):
\[(ad,b,c)=(a,bd,c)=(a,b,c)d.\]

We will use the following notation for a left annihilator of a subset $M$ in an algebra $A$:
\[Ann_l(M)=\{a\in A\,|\, aM=0\}.\]

\textbf{Lemma 3.} \textit{Let $I$ be an ideal in a Novikov algebra $A$, $n\in A$ and $(I,I,n)=0$. Then $(A,A,n)\subset\mathrm{Ann}_l T(I)$. If $n\in N(I)$, then $(A,n,A)=(n,A,A)\subset\mathrm{Ann}_l T(I)$.}

\textbf{Proof.} Indeed, let $a,b\in A$, $i,j\in I$. We have 
\[(a,b,n)((ij)j)=(a((ij)j),b,n)=-((a,ij,j),b,n)+((a(ij))j,b,n)\]
\[-((a,i,j),bj,n)=0.\]

It means that $(A,A,n)\subset\mathrm{Ann}_l T(I)$.

Let $n\in N(I)$. The next goal is to prove $(n,A,A)\subset\mathrm{Ann}_l T(I)$. First, we recall (\cite{P2022}) that $n(I,I,I)=0$. Then we have

\[(n,a,b)((ij)j)=((n(ij))j,a,b)= (nij,aj,b)=((n,aj,b)i)j\]
\[=(((n\cdot aj)b)i)j-((n(aj\cdot b))i)j\]
\[=(n(aj),b,i)j+((n\cdot aj)(bi))j-((n(aj\cdot b))i)j\]
\[=(n(aj)(bi))j-(n(aj\cdot b)i)j\]
\[=(n(aj\cdot bi))j-(n(aj\cdot b)i)j=(n(aj\cdot bi)-n(aj\cdot b)i)j\]
\[=(n(a\cdot bi)j)j-(n(ab\cdot j)i)j=(n(a\cdot bi-ab\cdot i)j)j=-(n(a,b,i)j)j\]
\[=-n(aj,bj,i)=0.\]

Thus, $(n,a,b)((ij)j)=0$ and $(n,A,A)\subset \mathrm{Ann}_l T(I)$. \qed

\medskip

Now we can strengthen one of the results of \cite{P2022}.\medskip

\textbf{Theorem 1.} \textit{Let $A$ be a prime nonassociative Novikov algebra with a nonzero ideal $I$. Then $I$ is a prime non-associative Novikov algebra.}

\textbf{Proof.} By Theorem 3 in \cite{P2022}, we obtain that $I$ is a prime Novikov algebra. Let $n\in N(I)$. By Lemma 2, $T(I)\neq 0$. By \cite{P2022} (Lemma 2), the left annihilator of a left ideal is an ideal. Then, by the primality of the algebra $A$, we obtain that $\mathrm{Ann}_lT(I)=0$. By Lemma 3 $(A,A,n)\subset\mathrm{Ann}_l T(I)=0$ and $(n,A,A)\subset\mathrm{Ann}_l T(I)=0 $. Thus $n\in N(A)=0$, so $N(I)=0$. In particular, the algebra $I$ is non-associative. \qed

\medskip

\textbf{Definition 4.} The intersection of all nonzero ideals of an algebra is called its \textbf{\textit{heart}}. An algebra is \textbf{\textit{subdirectly irreducible}} if its heart is nonzero. The heart $H$ of a subdirectly irreducible algebra is \textbf{\textit{idempotent}} if $H^2=H$. Since a subdirectly irreducible algebra with an idempotent heart is prime, the following corollary is true. Part of this corollary followed earlier from the results of \cite{ShestZhang}. \medskip

\textbf{Corollary 1.} \textit{Let $A$ be a subdirectly irreducible nonassociative Novikov algebra with an idempotent heart. Then the heart is a simple non-associative Novikov algebra.}

\medskip

\section{Andrunakievich radical}

\textbf{Definition 5.} A class of Novikov algebras $\mathcal{R}$ is called \textbf{\textit{radical}} if the following conditions are satisfied.

1) The homomorphic image of an algebra in $\mathcal{R}$ lies in $\mathcal{R}$.

2) Every Novikov algebra $A$ contains an ideal $\mathcal{R}(A)$ from $\mathcal{R}$, which contains all ideals in $A$ from $\mathcal{R}$.

3) The quotient algebra $A/\mathcal{R}(A)$ does not contain nonzero ideals from $\mathcal{R}$.\medskip

For completeness, we present the construction of the Andrunakievich radical from \cite{Zhevl}, restricting it to Novikov algebras.

Let $\mathcal{B}$ be the class of all subdirectly irreducible Novikov algebras with idempotent heart. An ideal $I$ of the Novikov algebra $A$ is called $\mathcal{B}$-ideal if the quotient algebra $A/I$ belongs to the class $\mathcal{B}$. The class $\mathcal{A}$ of Novikov algebras that do not map homomorphically onto algebras in the class $\mathcal{B}$ is radical (proved in \cite{Zhevl}). The largest ideal of the algebra $A$ contained in the class $\mathcal{A}$ is denoted by $\mathcal{A}(A)$ and is called the \textbf{\textit{Andrunakievich radical}} of the algebra $A$. \medskip

\textbf{Definition 6.} A radical class $\mathcal{R}$ is called \textbf{\textit{hereditary}} if for every ideal $I$ in $A$ we have $\mathcal{R}(I)=I\cap \mathcal{R}(A)$.


Let us prove some propositions. Their analogues for alternative algebras were proved in \cite{Zhevl} and our proofs are basically the same. 

\medskip\textbf{Proposition 1.} \textit{Let $A$ be a Novikov algebra, $I$ is an ideal in $A$ and $J$ is a $\mathcal{B}$-ideal in $I$. Then there exists a $\mathcal{B}$-ideal $K$ in $A$ with a condition $K\cap I=J$.}

\textbf{Proof.} An algebra $I/J$ is semiprime, so by Lemma~6 in \cite{P2022} $J$ is an ideal in $A$. We will use notations $\overline{A}=A/J$ and $\overline{I}=I/J$. The algebra $\overline{I}$ is a subdirectly irreducible Novikov algebra with idempotent heart. By Zorn's lemma we have an ideal $\overline{K}$ in $\overline{A}$, which is maximal with condition $\overline{I}\cap\overline{K}=\overline{0}$. It means that in an algebra $\overline{\overline{A}}=\overline{A}/\overline{K}$ all ideals have a nonzero intersection with an image of $\overline{I}$. Then the heart $H(\overline{\overline{A}})$ contains the heart of an image of $\overline{I}$ which is isomorphic to the heart of $\overline{I}$. But the heart of $\overline{I}$ is nonzero and idempotent. So, $\overline{\overline{A}}$ is in $\mathcal{B}$. 

Let $K$ be a preimage of an ideal $\overline{K}$ in $A$. Then we have $\overline{\overline{A}}=\overline{A}/\overline{K}\simeq A/K$. So, $K$ is a $\mathcal{B}$-ideal in $A$ and $K\cap I= J$. \qed

\medskip\textbf{Proposition 2.} \textit{Let $A$ be a Novikov algebra and $I$ be an ideal in $A$. Then $\mathcal{A}(I)$ is an ideal in $A$.}

\textbf{Proof.} We will use a notation $M=\mathcal{A}(I)$. 

a) Suppose that $M$ does not contain $Ma$ for some $a\in A$. Then $M\subsetneq Ma+M$ and $Ma+M$ is an ideal in $I$ by (\cite{P2022}, Lemma 5.a). Suppose that $(Ma)^2\subseteq M$. Let us define a map $\varphi:M\to (M+Ma)/M$, $\varphi(m)=ma+M$. If $m,n\in M$ then 
\[\varphi(m+n)=\varphi(m)+\varphi(n).\]
By right commutativity we have
\[\varphi(mn)=(mn)a+M=(ma)n+M=0+M.\]
We supposed that $(Ma)^2\subseteq M$, so
\[\varphi(m)\varphi(n)=(ma)(na)+M=0+M.\]
It means that $\varphi$ is a homomorphism of $M$ onto an ideal $(M+Ma)/M$ of an algebra $I/M$. It is a contradiction because $\mathcal{A}(I/M)=0$ and $M\in\mathcal{A}$. So $M$ does not contains $(Ma)^2$. 

Let us prove that $M+(Ma)^2$ is an ideal in $I$ and $M+(Ma)^2\subseteq M+Ma$. Indeed, $M+Ma$ is an ideal, so $M+(Ma)^2\subseteq M+Ma$. We have
\begin{multline*}(M+(Ma)^2))I\subseteq M+((Ma)(Ma))I\subseteq M+((Ma)I)(Ma)\subseteq \\ \subseteq M+(Ma)^2.
\end{multline*}
\begin{multline*}
    I(M+(Ma)^2)\subseteq M+I((Ma)(Ma)) \\ \subseteq M+(I,Ma,Ma)+(I(Ma))(Ma) \\ \subseteq M+(M,I,Ma)a+(Ma+M)(Ma) \\ \subseteq M+(Ma)^2+((MI)(Ma))a+(M(I(Ma)))a \\ \subseteq M+(Ma)^2+(Ma)^2I+Ma(Ma+M)\subseteq M+(Ma)^2.
\end{multline*}
So, $M+(Ma)^2$ is an ideal in $I$. 

Lemma 5.c in \cite{P2022} implies that $(M+(Ma)^2)/M$ is a trivial ideal in an algebra $I/M$. But $\mathcal{A}(I/M)=0$, so we have a contradiction. It means that $M$ contains $Ma$.

\medskip b) Suppose that $M$ does not contain $aM$ for some $a\in A$. Then $M\subsetneq aM+M$ and $aM+M$ is an ideal in $I$ by (\cite{P2022}, Lemma 5.a). Suppose that $(aM)^2\subseteq M$. Let us define a map $\varphi:M\to (M+aM)/M$, $\varphi(m)=am+M$. If $m,n\in M$ then 
\[\varphi(m+n)=\varphi(m)+\varphi(n).\]
We have
\[\varphi(mn)=a(mn)+M=(a,m,n)+(am)n+M=(m,a,n)+M=0+M.\]
We supposed that $(aM)^2\subseteq M$, so
\[\varphi(m)\varphi(n)=(am)(an)+M=0+M.\]
It means that $\varphi$ is a homomorphism of $M$ onto an ideal $(M+aM)/M$ of an algebra $I/M$. It is a contradiction because $\mathcal{A}(I/M)=0$ and $M\in\mathcal{A}$. So $M$ does not contain $(aM)^2$. 

Let us prove that $M+(aM)^2$ is an ideal in $I$ and $M+(aM)^2\subseteq M+aM$. Indeed, $M+aM$ is an ideal, so $M+(aM)^2\subseteq M+aM$. We have
\begin{multline*}(M+(aM)^2))I\subseteq M+((aM)(aM))I\subseteq M+((aM)I)(aM) \\ \subseteq M+(aM)^2.
\end{multline*}
\begin{multline*}
    I(M+(aM)^2)\subseteq M+I((aM)(aM)) \\ \subseteq M+(I,aM,aM)+(I(aM))(aM) \\ \subseteq M+(a,I,aM)M+(aM+M)(aM) \\ \subseteq M+(aM)^2+ ((aI)(aM))M+(a(I(aM)))M \\ \subseteq M+(aM)^2+((aM)(aM))I+(aM)(I(aM)) \\ \subseteq M+(aM)^2+(aM)M\subseteq M+(aM)^2\end{multline*}
So, $M+(aM)^2$ is an ideal in $I$. 

Lemma 5.e in \cite{P2022} implies that $(M+(aM)^2)/M$ is a trivial ideal in an algebra $I/M$. But $\mathcal{A}(I/M)=0$, so we have a contradiction. It means that $M$ contain $aM$. \qed

\medskip\textbf{Proposition 3.} \textit{In the class of Novikov algebras, the Andrunakievich radical is hereditary.}

\textbf{Proof.} Theorem 8.3 in \cite{Zhevl} states that a radical $\mathcal{R}$ is hereditary iff it satisfies two conditions:\\
(I) if $\mathcal{R}(A)=0$ and $I$ is an ideal in $A$ then $\mathcal{R}(I)=0$;\\
(II) if $\mathcal{R}(A)=A$ and $I$ is an ideal in $A$ then $\mathcal{R}(I)=I$.

Let us prove (I). Let $A$ be a Novikov algebra, $\mathcal{A}(A)=0$ and $I$ is a nonzero ideal in $A$. Then $\mathcal{A}(I)$ is an ideal in $A$ by Proposition~2. But $\mathcal{A}(\mathcal{A}(I))=\mathcal{A}(I)$, so $\mathcal{A}(I)\subseteq \mathcal{A}(A)=0$.

\medskip Let us prove (II). Let $A\in\mathcal{A}$ and $I$ is an ideal in $A$. If $I\notin\mathcal{A}$ then $I$ contains a nonzero $\mathcal{B}$-ideal $J\neq I$. By Proposition~1 we have a $\mathcal{B}$-ideal $K$ in $A$ such that $K\cap I=J$. It means that $A/K\in\mathcal{B}$ and $K\neq A$. We have a contradiction ($\mathcal{A}(A/K)=A/K)$. 

So, the Andrunakievich radical is hereditary in Novikov algebras. \qed

\medskip\textbf{Proposition 4.} \textit{Let $A$ be a Novikov algebra and $A\in\mathcal{B}$ and $I$ is a nonzero ideal in $A$. Then $I\in\mathcal{B}$ and $H(I)=H(A)$.}

\textbf{Proof.} An algebra $H(A)$ is simple by \cite{ShestZhang} and $H(A)\subseteq I$. Zorn's lemma allows us to choose a maximal ideal $K$ in $I$ with a condition $K\cap H(A)=0$. An algebra $I/K$ is semiprime, because each ideal in $I/K$ contains an image $\overline{H(A)}$ of $H(A)$. Then $K$ is an ideal in $A$ by (\cite{P2022}, Lemma 6), so $K=0$. It means that each nonzero ideal in $I$ has a nonzero intersection with $H(A)$. But $H(A)$ is simple, so each nonzero ideal in $I$ contains $H(A)$. \qed

Now we are ready to prove the following statement.

\medskip\textbf{Theorem 2.} \textit{The Andrunakievich radical $\mathcal{A}(A)$ of the Novikov algebra $A$ is equal to the intersection of all its $\mathcal{B}$-ideals, and the quotient algebra $A/\mathcal{A}(A)$ is a subdirect product of subdirectly irreducible algebras with an idempotent heart.}

\textbf{Proof.} If $C$ is an algebra in $\mathcal{B}$ and $L$ is an ideal in $C$ then $L$ is in $\mathcal{B}$ by Proposition~4. It means that $\mathcal{A}(L)\neq L$ and $\mathcal{A}(C)=0$.

Let $J$ be a $\mathcal{B}$-ideal in $A$. By the statement above we have $\mathcal{A}(A/J)=0$. It means that $\mathcal{A}(A)\subseteq J$.

Let $I$ be an intersection of all $\mathcal{B}$-ideals in an algebra $A$.  Then we have $\mathcal{A}(A)\subseteq I$.

If $I\neq \mathcal{A}(I)$ then $I$ contains a nonzero $\mathcal{B}$-ideal $J\neq I$. By Proposition~1 there exists a $\mathcal{B}$-ideal $K$ in $A$, such that $K\cap I=J\neq I$. We have a contradiction, $I$ is contained in all $\mathcal{B}$-ideals in $A$. So, $I\subseteq \mathcal{A}(A)$. Last statement is straightforward. \qed

\medskip

Since the sum of ideals is an ideal, then in finite dimensional algebras one can define a solvable radical as the largest solvable ideal. In the variety of Novikov algebras this ideal coincides with the Baer radical $B(A)$ whose existence was proved in \cite{P2022}. 




 Now we can prove that the Baer and Andrunakievich radicals coincide in the finite dimensional case.

\medskip\textbf{Theorem 3.} \textit{In a finite dimensional Novikov algebra, the Baer radical coincides with the Andrunakievich radical.}

\textbf{Proof.} Let $B$ be a prime finite dimensional Novikov algebra. If $H(B)=0$ then we have a finite set of nonzero ideals in $B$, which intersection is zero. It means that there exist nonzero ideals $I$ and $J$ in $B$, such that $I\cap J=0$. It is a contradiction with primality, so $H(B)\neq 0$ and $B\in\mathcal{B}$. 

Let $A$ be a Novikov algebra. By Theorem~2 $A/\mathcal{A}(A)$ is a subdirect product of prime Novikov algebras, so $A/\mathcal{A}(A)$ is a semiprime algebra and $B(A)\subseteq \mathcal{A}(A)$.

Consider the algebra $C=A/B(A)$. If $\mathcal{A}(A)$ does not lie in $B(A)$, then the algebra $C$ has a nonzero Andrunakievich radical $\mathcal{A}(C)\neq 0$. The algebra $C$ is a subdirect product of prime algebras $C/C_i$, where $C_i$ are ideals in $C$, $i\in \{1,\dots,n\}$. But $\mathcal{A}(C)$ has a nonzero image in some $C/C_k$, so $\overline{C}=C/C_k$ is an algebra in $\mathcal{B}$ and it has a nonzero Andrunakievich radical $\mathcal{A}(\overline{C})$. It means that an algebra $\mathcal{A}(\overline{C})$ is an algebra in $\mathcal{B}$ (by Proposition~4) and in $\mathcal{A}$. It is a contradiction, so $\mathcal{A}(A)\subseteq B(A)$. \qed 



\medskip

Note that in associative and commutative algebras the Andrunakievich radical coincides with the Jacobson radical, so that in infinite-dimensional Novikov algebras the Andrunakievich and Baer radicals differ.

\medskip

    \section{Quasiregular radical}

\textbf{Definition 7.} An element $x\in A$ is called \textit{\textbf{left quasiregular}} if there exists $y\in A$ such that $x+y=yx$. An element $x\in A$ is called \textit{\textbf{right quasiregular}} if there exists $y\in A$ such that $x+y=xy$. An element is called \textit{\textbf{quasiregular}} if there exists $y\in A$ such that $x+y=xy=yx$. An algebra is called \textit{\textbf{quasiregular}} (\textit{\textbf{right quasiregular}}, \textit{\textbf{left quasiregular}}) if all its elements are quasiregular (right quasiregular, left quasiregular).\medskip

\textbf{Example 1.} \textit{Consider the two-dimensional Novikov algebra $A=Fa+Fb$ with multiplication $a^2=b^2=ba=0$, $ab=b$, \cite{Zelm}. It is easy to see that the elements $\alpha a$ and $\beta b$ are quasiregular for every $\alpha,\beta\in F$, that is, the spaces $Fa$ and $Fb$ are quasiregular. However, it is directly verified that the element $a+b$ is not right quasiregular. Thus, the largest quasiregular (right quasiregular) ideal in $A$ is $Fb$. In particular, a solvable radical need not be quasiregular or right quasiregular in finite dimensional Novikov algebras. Since an algebra with zero multiplication is obviously quasiregular, the property of quasiregularity (right quasiregularity) is not radical in the class of finite dimensional Novikov algebras, namely, the following property does not hold: if $I$ and $A/I$ are radical, then $A$ is radical.}

\medskip

\textbf{Lemma 4.} \textit{The right-nilpotent algebra is left quasiregular. A left-nilpotent algebra is right quasiregular.}

\textbf{Proof.} Indeed, let $A$ be a right-nilpotent algebra, $x\in A$. Define $x^m=x^{m-1}\cdot x$ and let $x^n=0$. Then
\[-x+(x-x^2+\dots+(-1)^{n}x^{n-1})=(x-x^2+\dots+(-1)^{n}x^{n-1} )\cdot (-x).\]
Thus the element $-x$ is left quasiregular. But this is true for every $x\in A$. In other words, $y$ is left quasiregular for every $y\in A$. Similarly, it can be shown that every element of a left-nilpotent algebra is right quasiregular. \qed

\medskip

I.~P.~Shestakov and Z.~Zhang proved \cite{ShestZhang} that the right-nilpotency of the Novikov algebra is equivalent to its solvability (and is equivalent to the nilpotency of its square). Thus, Lemma 4 implies that the Baer radical of a finite dimensional Novikov algebra is left quasiregular. In fact, under some constraint on the field, the converse is also true.

\medskip

\textbf{Example 2.} \textit{Let $F$ be a field of characteristic $p>2$. Consider the algebra $A_{p^n}(a,b)$ with a basis $\{y_{-1},y_0,\dots,y_{p^n-2}\}$ and the following multiplication:
\[y_{-1}y_{-1}=ay_{p^n-2},\]
\[y_{-1}y_0=y_{-1}+by_{p^n-2},\]
\[y_iy_j=C_{i+j+1}^jy_{i+j},\]
where $a,b$ are fixed elements from $F$ and $y_k=0$ for $k>p^n-2$.}

\medskip In \cite{Xu} it is proved that over an algebraically closed field of characteristic $p>2$ every simple Novikov algebra is isomorphic to some $A_{p^n}(a,b)$.

\medskip

\textbf{Theorem 4.} \textit{Let $F$ be a field of characteristic 0 or an algebraically closed field of characteristic $p>2$. A finite dimensional Novikov algebra over the field $F$ is left quasiregular if and only if it is solvable.}

\textbf{Proof.} In one direction, the statement is true by Lemma~4. Let $A$ be a left quasiregular finite dimensional algebra. We prove by induction on dimension that $A$ is solvable. If $A=Fx$ and $x^2=\beta x\neq 0$ then $(\frac{1}{\beta}x)^2=\frac{1}{\beta^2}\beta x = \frac{1}{\beta}x$. It means that $e=\frac{x}{\beta}$ is an idempotent, $e^2=e$ and $A=Fe$. There exists $\alpha\in F$, such that $e+\alpha e = e\cdot \alpha e=\alpha e$, so $e=0$, it is a contradiction. So $x^2=0$ and $A$ is trivial.

Suppose $A$ has its proper non-zero ideal $I$. Then, by the inductive hypothesis, the algebras $I$ and $A/I$ are solvable, so that the algebra $A$ is solvable. Thus, we can assume that the algebra $A$ is simple. In \cite{Zelm} (see also \cite{DIU}) it is proved that over a field of characteristic 0 the simple finite dimensional Novikov algebra is a field. It is easy to see that an identity 1 is not left quasiregular, so a field is not left quasiregular. If $F$ is an algebraically closed field of characteristic $p>2$, then $A\simeq A_{p^n}(a,b)$ for some $n\in\mathbb{N}$, $a, b\in A$. Suppose $A$ is left quasiregular. Then the element $y_0$ is left quasiregular, that is, $y_0+x=xy_0$ for some $x\in A$. Let $x=\sum\limits_{j=-1}^{p^n-2} \alpha_j y_j$. Then
\[y_0+x=\sum\limits_{j=-1}^{p^{n}-2}\alpha_j y_j y_0=x+\alpha_{-1}by_{p^n-2},\]
whence it follows that $y_0=\alpha_{-1}by_{p^n-2}$, which is impossible. Contradiction. \qed

The question remains whether there exists the largest left quasiregular ideal in infinite dimensional Novikov algebras and whether it is radical in the class of Novikov algebras.

    \section{Acknowledgements}

The author expresses his gratitude to P.~S.~Kolesnikov, whose remarks helped to improve this paper. The author is grateful to the referee, whose remarks helped to improve this paper and to fix some mistakes.

\end{document}